\documentclass{ifacconf}%
\usepackage{amsmath}
\usepackage{graphicx}
\usepackage{natbib}
\usepackage{amsfonts}
\usepackage{amssymb}%
\providecommand{\U}[1]{\protect\rule{.1in}{.1in}}
\makeatletter
\renewcommand{\@cite}[1]{#1}
\def\@biblabel#1{\hspace*{-\labelsep}}
\makeatother
\begin{document}
\begin{frontmatter}
\title{Infinitesimal Perturbation Analysis for Quasi-Dynamic Traffic Light Controllers\thanksref{footnoteinfo}}
\thanks[footnoteinfo]{The authors' work is supported in part by NSF under
Grant CNS-1139021, by AFOSR under grant FA9550-12-1-0113, by
ONR under grant N00014-09-1-1051, and by ARO under Grant W911NF-11-1-0227.}
\author[First]{Julia L. Fleck}
\author[First]{Christos G. Cassandras}
\address[First]{Division of Systems Engineering and Center for Information and Systems Engineering,
Boston University, Brookline, MA 02446 USA (e-mail: jfleck@bu.edu, cgc@bu.edu)}
\begin{abstract}                
We consider the traffic light control problem for a single intersection modeled as a stochastic hybrid system. We
study a quasi-dynamic policy based on partial state information defined by detecting whether vehicle backlogs
are above or below certain controllable thresholds. Using Infinitesimal Perturbation Analysis (IPA), we derive online
gradient estimators of a cost metric with respect to these threshold parameters and use these estimators
to iteratively adjust the threshold values through a standard gradient-based algorithm so as to improve overall system
performance under various traffic conditions. Results obtained by applying this methodology to a simulated urban
setting are also included.
\end{abstract}
\end{frontmatter}

\section{Introduction}

The Traffic Light Control (TLC) problem consists of adjusting green and red
light cycles in order to control the traffic flow through an intersection and,
more generally, through a set of intersections and traffic lights. The
ultimate objective is to minimize congestion (hence delays experienced by
drivers) at a particular intersection, as well as an entire area consisting of
multiple intersections with traffic lights. Recent technological developments
have made it possible to collect and process traffic data so that they may be
applied in solving the TLC problem in real time. Fundamentally, TLC is a form
of scheduling for systems operating through simple switching control actions.
Numerous solution algorithms have been proposed and we briefly review some of
them. \cite{Porche1996} used a decision tree model with a Rolling Horizon
Dynamic Programming (RHDP) approach, while \cite{Dujardin2011} proposed a
multiobjective Mixed Integer Linear Programming (MILP) formulation. Optimal
TLC was also stated as a special case of an Extended Linear Complementarity
Problem (ELCP) by \cite{DeSchutter1999}, and formulated as a hybrid system
optimization problem by \cite{Zhao2003}. \cite{Yu2006} modeled a traffic light
intersection as a Markov Decision Process (MDP) and a game theoretic approach
was applied to a finite controlled Markov chain model by \cite{Alvarez2010}.
Relying on sensor information regarding traffic congestion, \cite{Choi2002}
developed a first-order Sugeno fuzzy model and incorporated it into a fuzzy
logic controller. Perturbation analysis techniques were used by
\cite{Head1996} and \cite{Fu2003} for modeling a traffic light intersection as
a stochastic Discrete Event System (DES), while an Infinitesimal Perturbation
Analysis (IPA) approach, using a Stochastic Flow Model (SFM) to represent the
queue content dynamics of roads at an intersection, was presented in
[\cite{Panayiotou2005}].

Our work is also based on modeling traffic flow through an intersection
controlled by switching traffic lights as an SFM, which conveniently captures
the system's inherent hybrid nature: while traffic light switches exhibit
event-driven dynamics, the flow of vehicles through an intersection is best
represented through time-driven dynamics. In [\cite{Geng2012}], IPA was
applied with respect to controllable green and red cycle lengths for a single
isolated intersection and in [\cite{Geng2013a}] for multiple intersections.
Traffic flow rates need not be restricted to take on deterministic values, but
may be treated as stochastic processes (see [\cite{Cassandras2002}]), which
are suited to represent the continuous and random variations in traffic
conditions. Using the general IPA theory for Stochastic Hybrid Systems (SHS)
in [\cite{Wardi2010}] and [\cite{Cassandras2010}], on-line gradients of
performance measures are estimated with respect to several controllable
parameters with only minor technical conditions imposed on the random
processes that define input and output flows. These IPA estimates have been
shown to be unbiased, even in the presence of blocking due to limited resource
capacities and of feedback control (see [\cite{Yao2011}]). It should be
emphasized that IPA is not used to estimate performance measures, but only
their gradients, which may be subsequently incorporated into standard
gradient-based algorithms in order to effectively control parameters of interest.

In contrast to earlier work where the adjustment of light cycles did not make
use of real-time state information, \cite{Geng2013b} proposed a quasi-dynamic
control setting in which partial state information is used conditioned upon a
given queue content threshold being reached. In this paper, we draw upon this
setting, but rather than controlling the light cycle lengths as in
[\cite{Geng2013b}], here we focus on the threshold parameters and derive IPA
performance measure estimators necessary to optimize these parameters, while
assuming fixed cycle lengths. Our goal is to compare the relative effects of
the threshold parameters and the light cycle length parameters on our
objective function, build upon these results, and ultimately control both the
light cycle lengths and the queue content thresholds simultaneously.

In Section 2, we formulate the TLC\ problem for a single intersection and
present the modeling framework used throughout our analysis for controlling
vehicle queue thresholds. Section 3 details the derivation of an
IPA\ estimator for the cost function gradient with respect to a controllable
parameter vector defined by these thresholds. The IPA estimator is then
incorporated into a gradient-based optimization algorithm and we include
simulation results in Section 4, showing how the proposed quasi-dynamic
control offers considerable improvement over prior results.

\section{Problem Formulation}

\begin{figure}
\begin{center}%
$\includegraphics{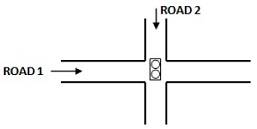}$\caption{A single traffic
light intersection with two cross-roads}\label{(fig): Single Intersection}%
\end{center}
\end{figure}

The system we consider comprises a single intersection, as shown in Fig.
\ref{(fig): Single Intersection}. For simplicity, left-turn and right-turn
traffic flows are not considered and yellow light cycles are implicitly
accounted for within a red light cycle. We assign to each queue $i$ a
guaranteed minimum GREEN cycle length $\theta_{i,\min}$, and a maximum length
$\theta_{i,\max}$ which (in contrast to [\cite{Geng2013b}]) we assume to be
fixed. We define a state vector $x(t)=[x_{1}(t),x_{2}(t)]$ where $x_{i}%
(t)\in\mathbb{R}^{+}$ is the content of queue $i$. For each queue $i$, we also
define a \textquotedblleft clock\textquotedblright\ state variable $z_{i}(t)$,
$i=1,2$, which measures the time since the last switch from RED to GREEN of
the traffic light for queue $i$, so that $z_{i}(t)\in\lbrack0,\theta_{i,\max
}]$. Setting $z(t)=[z_{1}(t),z_{2}(t)]$, the complete system state vector is
$[x(t),z(t)]$. Within the \emph{quasi-dynamic} setting considered in this
work, the controllable parameter vector of interest is given by $\mathbf{s}%
=\left[  s_{1},s_{2}\right]  $, where $s_{n}\in\Re^{+}$ is a queue content
threshold value for road $n=1,2$. The notation $x(\mathbf{s},t)=[x_{1}%
(\mathbf{s},t),x_{2}(\mathbf{s},t)]$ is used to stress the dependence of the
state on these threshold parameters. However, for notational simplicity, we
will henceforth write $x(t)$ when no confusion arises; the same applies to
$z(t)$.

Let us now partition the queue content state space into the following four
regions:
\[
X_{0}=\{(x_{1},x_{2}):x_{1}(t)<s_{1},\text{ }x_{2}(t)<s_{2}\}
\]%
\[
X_{1}=\{(x_{1},x_{2}):x_{1}(t)<s_{1},\text{ }x_{2}(t)\geq s_{2}\}
\]%
\[
X_{2}=\{(x_{1},x_{2}):x_{1}(t)\geq s_{1},\text{ }x_{2}(t)<s_{2}\}
\]%
\[
X_{3}=\{(x_{1},x_{2}):x_{1}(t)\geq s_{1},\text{ }x_{2}(t)\geq s_{2}\}
\]
%

At any time $t$, the feasible control set for the traffic light controller is
$U=\{1,2\}$ and the control is defined as:%
\begin{equation}
u\left(  x(t),z(t)\right)  \equiv\left\{
\begin{array}
[c]{c}%
1\\
2
\end{array}
\right.
\begin{array}
[c]{l}%
\text{i.e., set road }1\text{ GREEN, road }2\text{ RED}\\
\text{i.e., set road }2\text{ GREEN, road }1\text{ RED}%
\end{array}
\label{(eq): Control definition}%
\end{equation}
A \emph{dynamic} controller is one that makes full use of the state
information $z(t)$ and $x(t)$. Obviously, $z(t)$ is the controller's known
internal state, but the queue content state is generally not observable. We
assume, however, that it is \emph{partially} observable. Specifically, we can
only observe whether $x_{i}(t)$ is below or above some threshold $s_{i},i=1,2$
(this is consistent with actual traffic systems where sensors (typically,
inductive loop detectors) are installed at each road near the intersection).
In this context, we shall define a \emph{quasi-dynamic} controller of the form
$u\left(  X(t),z(t)\right)  $, with $X(t)\in\left\{  X_{0},X_{1},X_{2}%
,X_{3}\right\}  $, as follows:

For $X(t)\in\left\{  X_{0},X_{3}\right\}  $:%
\begin{equation}
u\left(  z(t)\right)  =\left\{
\begin{array}
[c]{c}%
1\\
2
\end{array}
\right.
\begin{array}
[c]{l}%
\text{if }z_{1}(t)\in\left(  0,\theta_{1,\max}\right)  \text{ and }%
z_{2}(t)=0\\
\text{otherwise}%
\end{array}
\label{(eq): Control policy I}%
\end{equation}

For $X(t)=X_{1}$:%
\begin{equation}
u\left(  z(t)\right)  =\left\{
\begin{array}
[c]{c}%
1\\
2
\end{array}
\right.
\begin{array}
[c]{l}%
\text{if }z_{1}(t)\in\left(  0,\theta_{1,\min}\right)  \text{ and }%
z_{2}(t)=0\\
\text{otherwise}%
\end{array}
\label{(eq): Control policy II}%
\end{equation}

For $X(t)=X_{2}$:%
\begin{equation}
u\left(  z(t)\right)  =\left\{
\begin{array}
[c]{c}%
2\\
1
\end{array}
\right.
\begin{array}
[c]{l}%
\text{if }z_{2}(t)\in\left(  0,\theta_{2,\min}\right)  \text{ and }%
z_{1}(t)=0\\
\text{otherwise}%
\end{array}
\label{(eq): Control policy III}%
\end{equation}
This is a simple form of hysteresis control to ensure that the $i$th traffic
flow always receives a minimum GREEN light cycle $\theta_{i,\min}$. Clearly,
the GREEN light cycle may be dynamically interrupted anytime after
$\theta_{i,\min}$ based on the partial state feedback provided through $X(t)$.
For instance, if a transition into $X_{1}$ occurs while $u\left(
X(t),z(t)\right)  =1$ and $z_{1}(t)>\theta_{1,\min}$, then the light switches
from GREEN to RED for road 1 in order to accommodate an increasing backlog at
road 2. For notational simplicity, we will write $u(t)$ when no confusion
arises, as we do with $x(t)$, $z(t).$

The stochastic processes involved in this system are defined on a common
probability space $\left(  \Omega,F,P\right)  $. The arrival flow processes
are $\{\alpha_{n}(t)\}$, $n=1,2$, where $\alpha_{n}(t)$ is the instantaneous
vehicle arrival rate at time $t$. The departure flow process on road $n$ is
defined as:%
\begin{equation}
\beta_{n}(t)=\left\{
\begin{array}
[c]{l}%
h_{n}(X(t),z(t),t)\\
\alpha_{n}(t)\\
0
\end{array}
\right.
\begin{array}
[c]{l}%
\text{if }x_{n}(t)>0\text{ and }u(t)=n\\
\text{if }x_{n}(t)=0\text{ and }u(t)=n\\
\text{otherwise}%
\end{array}
\label{(eq): Departure process}%
\end{equation}
where $h_{n}(X(t),z(t),t)$ is the instantaneous vehicle departure rate at time
$t$; for notational simplicity, we will write $h_{n}(t)$ when no confusion
arises. We can now write the dynamics of the state variables $x_{n}(t)$ and
$z_{n}(t)$ as follows, where we adopt the notation $\overline{n}$ to denote
the index of the road perpendicular to road $n=1,2$, and note that the symbols
$t^{+}$ ($t^{-}$, respectively) denote the time instant immediately following
(preceding, respectively) time $t$:%
\begin{equation}
\overset{\cdot}{x}_{n}(t)=\left\{
\begin{array}
[c]{l}%
\alpha_{n}(t)\\
0\\
\alpha_{n}(t)-h_{n}(t)
\end{array}
\right.
\begin{array}
[c]{l}%
\text{if }z_{n}(t)=0\\
\text{if }x_{n}(t)=0\text{ and }\alpha_{n}(t)\leq h_{n}(t)\\
\text{otherwise}%
\end{array}
\label{(eq): State dynamics}%
\end{equation}%
\begin{equation}
\overset{\cdot}{z}_{n}(t)=\left\{
\begin{array}
[c]{c}%
1\\
0
\end{array}
\right.
\begin{array}
[c]{l}%
\text{if }z_{\overline{n}}(t)=0\\
\text{otherwise}%
\end{array}
\label{(eq): Clock dynamics}%
\end{equation}

\[%
\begin{array}
[t]{l}%
z_{n}(t^{+})=0\\
\text{if }z_{n}(t)=\theta_{n,\max}\\
\text{or }z_{n}(t)=\theta_{n,\min},\text{ }x_{n}(t)<s_{n},\text{ }%
x_{\overline{n}}(t)\geq s_{\overline{n}}\\
\text{or }z_{n}(t)>\theta_{n,\min},\text{ }x_{n}(t^{-})>s_{n},x_{n}%
(t^{+})=s_{n},x_{\overline{n}}(t)\geq s_{\overline{n}}\\
\text{or }z_{n}(t)>\theta_{n,\min},\text{ }x_{n}(t)<s_{n},x_{\overline{n}%
}(t^{-})<s_{\overline{n}},x_{\overline{n}}(t^{+})=s_{\overline{n}}%
\end{array}
\]

In this context, the traffic light intersection in Fig.
\ref{(fig): Single Intersection} can be viewed as a hybrid system in which the
time-driven dynamics are given by (\ref{(eq): State dynamics}%
)-(\ref{(eq): Clock dynamics}) and the event-driven dynamics are associated
with light switches and with events that cause the value of $x_{n}(t)$ to
change from strictly positive to zero or vice-versa. It is then possible to
derive a Stochastic Hybrid Automaton (SHA) model as in [\cite{Geng2013b}]
containing 14 modes, which are defined by combinations of $x_{n}(t)$ and
$z_{n}(t)$ values. The event set for this SHA is $\Phi_{n}=\left\{
e_{1},e_{2},e_{3},e_{4},e_{5},e_{6},e_{7}\right\}  $, where $e_{1}$ is the
guard condition $\left[  x_{n}=s_{n}\text{ from below}\right]  $; $e_{2}$ is
the guard condition $\left[  x_{n}=s_{n}\text{ from above}\right]  $; $e_{3}$
is the guard condition $\left[  z_{n}=\theta_{n,\min}\right]  $; $e_{4}$ is
the guard condition $\left[  z_{n}=\theta_{n,\max}\right]  $; $e_{5}$ is the
guard condition $\left[  x_{n}=0\text{ from above}\right]  $; $e_{6}$ is a
switch in the sign of $\alpha_{n}(t)-h_{n}(t)$ from non-positive to strictly
positive; $e_{7}$ is a switch in the sign of $\alpha_{n}(t)$ from 0 to
strictly positive. Note that $e_{1},\ldots,e_{4}$ are the events that induce
light switches and, for easier reference, we rename them as $\zeta_{n}$,
$\gamma_{n}$, $\lambda_{n}$, and $\mu_{n}$, respectively, where the subscript
$n$ refers to the road where the event occurred. If we label light switching
events from RED to GREEN and GREEN to RED as $R2G_{n}$ and $G2R_{n}$,
respectively, we can specify the following rules for our hysteresis-based controller:

\begin{description}
\item[Rule 1] The occurrence of event $\zeta_{n}$, while $z_{\overline{n}%
}>\theta_{\bar{n},\min}$ and $x_{\overline{n}}<s_{\overline{n}}$, results in
event $R2G_{n}$.

\item[Rule 2] The occurrence of event $\gamma_{n}$, while $z_{n}%
>\theta_{n,\min}$ and $x_{\overline{n}}\geq s_{\overline{n}}$, results in
event $G2R_{n}$.

\item[Rule 3] The occurrence of event $\lambda_{n}$, while $x_{n}<s_{n}$ and
$x_{\overline{n}}\geq s_{\overline{n}}$, results in event $G2R_{n}$.

\item[Rule 4] The occurrence of event $\mu_{n}$ always results in event
$G2R_{n}$.
\end{description}

\begin{figure}
\begin{center}%
$\includegraphics{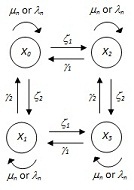}$%
\caption{Stochastic Hybrid Automaton for aggregate states $X(t)$ under quasi-dynamic control}\label{(fig): SHA}%
\end{center}
\end{figure}

\begin{figure}
\begin{center}%
$\includegraphics[width=8.4cm]{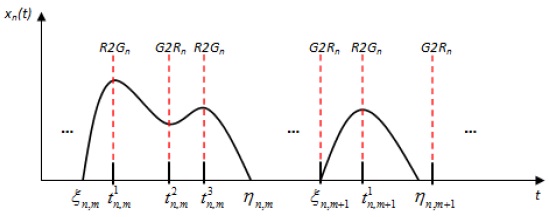}$\caption{Typical sample path of a
traffic light queue}\label{(fig): Sample Path}%
\end{center}
\end{figure}

A partial state transition diagram defined in terms of the aggregate queue
content states $X(t)$ is shown in Fig. \ref{(fig): SHA}. A complete state
transition diagram for this SHA is too complicated to draw and is not
necessary for IPA, which focuses on analyzing a typical sample path and
observable events in it, as shown in Fig. \ref{(fig): Sample Path}. Observe
that any such sample path consists of alternating Non-Empty Periods (NEPs) and
Empty Periods (EPs), which correspond to time intervals when $x_{n}(t)>0$
(i.e., queue $n$ is non-empty) and $x_{n}(t)=0$ (i.e., queue $n$ is empty),
respectively. Let us then label the events corresponding to the end and to the
start of an NEP as $E_{n}$ and $S_{n}$, respectively, and note that $E_{n}$ is
induced by event $e_{5}$, while $S_{n}$ may be induced by events $e_{6}$ or
$e_{7}$ or $G2R_{n}$.%

Let us denote the $m$th NEP in a sample path of queue $n$, by $\left[
\xi_{n,m},\eta_{n,m}\right)  $, where $\xi_{n,m}$, $m=1,2,\ldots$, is the time
of occurrence of the $m$th $S_{n}$ event and $\eta_{n,m}$ is the time of
occurrence of the $m$th $E_{n}$ event, as illustrated in Fig.
\ref{(fig): Sample Path}. Additionally, let the time of a light switching
event (either $R2G_{n}$ or $G2R_{n}$) within the $m$th NEP be denoted by
$t_{n,m}^{j}$, $j=1,...,J_{m}$.%

Recall that the purpose of our analysis is to apply IPA to sample path data in
order to obtain unbiased gradient estimates of a system performance measure
with respect to the controllable parameter vector\textbf{\ }$\mathbf{s}$ and
subsequently incorporate such estimates into a gradient-based optimization
scheme. In particular, we define a sample function which measures a weighted
mean of the queue lengths over a fixed time interval $[0,T]$:%
\begin{equation}
L\left(  \mathbf{s};x(0),z(0),T\right)  =\frac{1}{T}\underset{n=1}{\overset
{2}{\sum}}%
{\displaystyle\int\limits_{0}^{T}}
w_{n}x_{n}\left(  \mathbf{s},t\right)  dt\label{(eq): Sample function}%
\end{equation}
where $w_{n}$ is a weight associated with road $n$, and $x(0)$, $z(0)$ are
given initial conditions. Since $x_{n}(t)=0$ during EPs of road $n$%
,(\ref{(eq): Sample function}) can be rewritten as%
\begin{equation}
L\left(  \mathbf{s};x(0),z(0),T\right)  =\frac{1}{T}\overset{2}{\underset
{n=1}{\sum}}\underset{m=1}{\overset{M_{n}}{\sum}}%
{\displaystyle\int\limits_{\xi_{n,m}}^{\eta_{n,m}}}
w_{n}x_{n}\left(  \mathbf{s},t\right)
dt\label{(eq): Sample function simplified}%
\end{equation}
where $M_{n}$ is the total number of NEPs during the sample path of road $n$.
Finally, using $E$ to denote the usual expectation operator, let us define the
overall performance metric as%
\begin{equation}
J\left(  \mathbf{s};x(0),z(0),T\right)  =E\left[  L\left(  \mathbf{s}%
;x(0),z(0),T\right)  \right]  \label{(eq): Performance metric}%
\end{equation}
We note that it is not possible to derive a closed-form expression of
$J\left(  \mathbf{s};x(0),z(0),T\right)  $ without full knowledge of the
processes $\{\alpha_{n}(t)\}$ and $\{\beta_{n}(t)\}$. On the other hand, by
assuming only that $\alpha_{n}(t)$ and $\beta_{n}(t)$ are piecewise continuous
w.p. 1, we can successfully apply the IPA methodology developed for general
SHS by \cite{Cassandras2010} and obtain an estimate of $\nabla J\left(
\mathbf{s}\right)  $ by evaluating the sample gradient $\nabla L\left(
\mathbf{s}\right)  $. As we will see, no explicit knowledge of $\alpha_{n}%
(t)$, $\beta_{n}(t)$ is necessary to estimate $\nabla J\left(  \mathbf{s}%
\right)  $, which can then be used to improve current operating conditions or
(under certain conditions) to compute an optimal $\mathbf{s}^{\ast}$ through
an iterative optimization algorithm of the form%
\begin{equation}
s_{i,l+1}=s_{i,l}-\rho_{l}H_{i,l}\left(  \mathbf{s}_{l},x(0),T,\omega
_{l}\right)  \label{(eq): Optimization scheme}%
\end{equation}
where $\rho_{l}$ is the step size at the $l$th iteration, $l=0,1,...$, and
$\omega_{l}$ denotes a sample path from which data are extracted and used to
compute $H_{i,l}\left(  \mathbf{s}_{l},x(0),T,\omega_{l}\right)  $, which is
an estimate of $dJ/ds_{i}$. We will further assume that that the derivatives
$dL/ds_{i}$ exist w.p. 1 for all $s_{i}\in\Re^{+}$. It is also easy to check
that $L\left(  \mathbf{s}\right)  $ is Lipschitz continuous for $s_{i}\in
\Re^{+}$. Under these conditions, it has been shown by \cite{Cassandras2010}
that $dL/ds_{i}$ is an unbiased estimator of $dJ/ds_{i}$, $i=1,2$.

\section{Infinitesimal Perturbation Analysis}

For the sake of completeness we begin with a brief overview of the generalized
IPA framework developed for SHS in [\cite{Cassandras2010}]). Consider a sample
path of the system over $\left[  0,T\right]  $ and denote the time of
occurrence of the $k$th event (of any type) by $\tau_{k}\left(  \theta\right)
$, where $\theta$ is a scalar (for simplicity) controllable parameter of
interest. We shall also denote the state and event time derivatives with
respect to parameter $\theta$ as $x^{\prime}(\theta,t)\equiv\frac{\partial
x(\theta,t)}{\partial\theta}$ and $\tau_{k}^{\prime}(\theta)\equiv
\frac{\partial\tau_{k}\left(  \theta\right)  }{\partial\theta}$, respectively,
for $k=1,...,N$. The dynamics of $x(\theta,t)$ are fixed over any interevent
interval $[\tau_{k}(\theta),\tau_{k+1}(\theta))$ and we write $\dot{x}%
(\theta,t)=f_{k}(\theta,x,t)$ to represent the state dynamics over this
interval. Although we include $\theta$ as an argument in the expressions above
to stress dependence on the controllable parameter, we will subsequently drop
this for ease of notation as long as no confusion arises. It is shown in
[\cite{Cassandras2010}] that the state derivative satisifies%
\begin{equation}
\frac{d}{dt}x^{\prime}(t)=\frac{\partial f_{k}(t)}{\partial x}x^{\prime
}(t)+\frac{\partial f_{k}(t)}{\partial\theta}
\label{(eq): General state derivative}%
\end{equation}
with the following boundary condition:%
\begin{equation}
x^{\prime}(\tau_{k}^{+})=x^{\prime}(\tau_{k}^{-})+\left[  f_{k-1}(\tau_{k}%
^{-})-f_{k}(\tau_{k}^{+})\right]  \cdot\tau_{k}^{\prime}
\label{(eq): Boundary condition}%
\end{equation}
Knowledge of $\tau_{k}^{\prime}$ is, therefore, needed in order to evaluate
(\ref{(eq): Boundary condition}). Following the framework in
[\cite{Cassandras2010}], there are three types of events for a general
stochastic hybrid system. $(i)$ \emph{Exogenous Events.} These events cause a
discrete state transition independent of $\theta$ and satisfy $\tau
_{k}^{\prime}=0$. $(ii)$ \emph{Endogenous Events.} Such an event occurs at
time $\tau_{k}$ if there exists a continuously differentiable function
$g_{k}:\mathbb{R}^{n}\times\Theta\rightarrow\mathbb{R}$ such that $\tau
_{k}\ =\ \min\{t>\tau_{k-1}\ :\ g_{k}\left(  x\left(  \theta,t\right)
,\theta\right)  =0\}$, where the function $g_{k}$ normally corresponds to a
guard condition in a hybrid automaton. Taking derivatives with respect to
$\theta$, it is straightforward to obtain%
\begin{equation}
\tau_{k}^{\prime}=-\left[  \frac{\partial g_{k}}{\partial x}\cdot f_{k-1}%
(\tau_{k}^{-})\right]  ^{-1}\cdot\left(  \frac{\partial g_{k}}{\partial\theta
}+\frac{\partial g_{k}}{\partial x}\cdot x^{\prime}(\tau_{k}^{-})\right)
\label{(eq): Endogenous event}%
\end{equation}
where $\frac{\partial g_{k}}{\partial x}.f_{k-1}(\tau_{k}^{-})\neq0$. $(iii)$
\emph{Induced Events.} Such an event occurs at time $\tau_{k}$ if it is
triggered by the occurrence of another event at time $\tau_{m}\leq\tau_{k}$
(details can be found in [\cite{Cassandras2010}]).

Returning to our TLC problem, we define the derivatives of the states
$x_{n}(\mathbf{s,}t)$ and $z_{i}(\mathbf{s,}t)$ and event times $\tau
_{k}(\mathbf{s})$ with respect to $s_{i}$, $i=1,2$, as follows:
\begin{equation}
x_{n,i}^{\prime}(t)\equiv\frac{\partial x_{n}(\mathbf{s,}t)}{\partial s_{i}%
},\text{ }z_{i,i}^{\prime}(\mathbf{s,}t)\equiv\frac{\partial z_{i}%
(\mathbf{s,}t)}{\partial s_{i}},\text{ }\tau_{k,i}^{\prime}\equiv
\frac{\partial\tau_{k}(\mathbf{s})}{\partial s_{i}}\label{IPAnotation}%
\end{equation}
Observe that, based on (\ref{(eq): State dynamics}), $\frac{\partial
f_{n,k}(t)}{\partial x_{n}}=\frac{\partial f_{n,k}(t)}{\partial s_{i}}=0$,
$n,i=1,2,$ so that in (\ref{(eq): General state derivative}) we have $\frac
{d}{dt}x_{n}^{\prime}(t)=0$ for $t\in\lbrack\tau_{k},\tau_{k+1})$. Thus,
$x_{n}^{\prime}(t)=x_{n}^{\prime}(\tau_{k}^{+})$, $t\in\left[  \tau_{k}%
,\tau_{k+1}\right)  $. In what follows, we derive the IPA state and event time
derivatives for the events identified in our SHA model.

\subsection{State and Event Time Derivatives}

We shall proceed by considering each of the event types ($G2R_{n}$, $R2G_{n}$,
$E_{n}$, $S_{n}$) identified in the previous section and deriving the
corresponding event time and state derivatives. We begin with a general result
which applies to all light switching events $G2R_{n}$ and $R2G_{n}$. Let us
denote the time of the $j$th occurrence of a light switching event by
$\sigma_{j}$ and define its derivative with respect to the control parameters
as $\sigma_{j,i}^{\prime}\equiv\frac{\partial\sigma_{j}}{\partial s_{i}}$,
$i=1,2$.

\emph{Lemma 1:} The derivative $\sigma_{j,i}^{\prime}$, $i=1,2$, of light
switching event times $\sigma_{j}$, $j=1,2,\ldots$ with respect to the control
parameters $s_{1},s_{2}$ satisfies:%
\begin{equation}
\sigma_{j,i}^{\prime}=\left\{
\begin{array}
[c]{ll}%
\frac{\mathbf{1}\left[  n=i\right]  -x_{n,i}^{\prime}(\sigma_{j}^{-})}%
{\alpha_{n}(\sigma_{j})} & \text{if }\zeta_{n}\text{ occurs at }\sigma_{j}\\
\frac{\mathbf{1}\left[  n=i\right]  -x_{n,i}^{\prime}(\sigma_{j}^{-})}%
{\alpha_{n}(\sigma_{j})-h_{n}(\sigma_{j})} & \text{if }\gamma_{n}\text{ occurs
at }\sigma_{j}\\
\sigma_{j-1,i}^{\prime} & \text{otherwise}%
\end{array}
\right.  \label{lemma1}%
\end{equation}
where $\mathbf{1}\left[  \cdot\right]  $ is the usual indicator function.

\textbf{Proof}: We begin with a $G2R_{n}$ light switching event. This event is
induced by one of four possible endogenous events which we analyze separately
in what follows.

1. \emph{Event }$\zeta_{1}$\emph{ occurs at time }$\sigma_{j}$\emph{. }In this
case, a $G2R_{2}$ event occurs, hence also a $R2G_{1}$ event. Since road 1
must be undergoing a RED cycle within a NEP, it follows from
(\ref{(eq): Endogenous event}) with $g_{j}=x_{1}-s_{1}$ and
(\ref{(eq): State dynamics}) that $\sigma_{j,1}^{\prime}=\frac{1-x_{1,1}%
^{\prime}(\sigma_{j}^{-})}{\alpha_{1}(\sigma_{j})}$ and $\sigma_{j,2}^{\prime
}=\frac{-x_{1,2}^{\prime}(\sigma_{j}^{-})}{\alpha_{1}(\sigma_{j})}$.

2. \emph{Event }$\zeta_{2}$\emph{ occurs at time }$\sigma_{j}$\emph{. }This
results in a $G2R_{1}$ event and the same reasoning as above applies to verify
that $\sigma_{j,1}^{\prime}=\frac{-x_{2,1}^{\prime}(\sigma_{j}^{-})}%
{\alpha_{2}(\sigma_{j})}$ and $\sigma_{j,2}^{\prime}=\frac{1-x_{2,2}^{\prime
}(\sigma_{j}^{-})}{\alpha_{2}(\sigma_{j})}$.

3. \emph{Event }$\gamma_{1}$\emph{ occurs at time }$\sigma_{j}$\emph{. }This
results in a $G2R_{1}$ event. Moreover, since this a light switching event, it
follows from (\ref{(eq): Control policy II}) that $x_{1}(\sigma_{j}^{-}%
)>s_{1}$ and $x_{1}(\sigma_{j})=s_{1}$, which means that road 1 must be in a
NEP with $\beta_{1}(\sigma_{j})>0$. As a result, it follows from
(\ref{(eq): Endogenous event}) with $g_{j}=x_{1}-s_{1}$ and
(\ref{(eq): State dynamics}) that $\sigma_{j,1}^{\prime}=\frac{1-x_{1,1}%
^{\prime}(\sigma_{j}^{-})}{\alpha_{1}(\sigma_{j})-h_{1}(\sigma_{j})}$ and
$\sigma_{j,2}^{\prime}=\frac{-x_{1,2}^{\prime}(\sigma_{j}^{-})}{\alpha
_{1}(\sigma_{j})-h_{1}(\sigma_{j})}$.

4. \emph{Event }$\gamma_{2}$\emph{ occurs at time }$\sigma_{j}$\emph{. }This
results in a $G2R_{2}$ event and the same reasoning as above applies to verify
that $\sigma_{j,1}^{\prime}=\frac{-x_{2,1}^{\prime}(\sigma_{j}^{-})}%
{\alpha_{2}(\sigma_{j})-h_{2}(\sigma_{j})}$ and $\sigma_{j,2}^{\prime}%
=\frac{1-x_{1,2}^{\prime}(\sigma_{j}^{-})}{\alpha_{2}(\sigma_{j})-h_{2}%
(\sigma_{j})}$.

5. \emph{Event }$\lambda_{n}$, $n=1,2$,\emph{ occurs at time }$\sigma_{j}%
$\emph{. }Let $\Delta_{j}=\sigma_{j}-\sigma_{j-1}$, $j=1,2,\ldots$, where
(without loss of generality) we set $\sigma_{0}=0$. Therefore, we can write
$\sigma_{j}=\sigma_{j-1}+\Delta_{j}$, $j=1,2,\ldots$ Recall that, by
definition, whenever a light switching event is induced by $\lambda_{n}$ we
must have $\Delta_{j}=\theta_{n,\min}$, which is independent of $s_{1},s_{2}$.
Therefore, $\sigma_{j,i}^{\prime}=\sigma_{j-1,i}^{\prime}$ for all
$j=1,2,\ldots$ and $i=1,2$.

6. \emph{Event }$\mu_{n}$, $n=1,2$,\emph{ occurs at time }$\sigma_{j}$\emph{.
}This is similar to the previous case with $\Delta_{j}=\theta_{n,\max}$ and
once again we have $\sigma_{j,i}^{\prime}=\sigma_{j-1,i}^{\prime}$ for all
$j=1,2,\ldots$ and $i=1,2$.

This concludes the proof for a $G2R_{n}$ light switching event. The analysis
for a $R2G_{n}$ event is similar, due to the fact that the end of a RED cycle
on road $n$ ($R2G_{n}$ event) must coincide with the start of a RED cycle on
road $\overline{n}$ ($G2R_{\overline{n}}$ event). $\blacksquare$

We now proceed by considering each of the event types ($G2R_{n}$, $R2G_{n}$,
$E_{n}$, $S_{n}$).

\begin{description}
\item[\textbf{(1)}] Event $G2R_{n}$
\end{description}

Two cases must be considered: $(a)$ $G2R_{n}$ occurs at $\tau_{k}$ while road
$n$ is undergoing an NEP; $(b)$ $G2R_{n}$ occurs at $\tau_{k}$ while road $n$
is undergoing an EP. In case $(a)$, the fact that $x_{n}(\tau_{k}^{-})>0$
means that $f_{n,k-1}(\tau_{k}^{-})=\alpha_{n}(\tau_{k})-h_{n}(\tau_{k})$.
Additionally, since road $n$ is undergoing a RED cycle at time $\tau_{k}^{+}$,
we must have $f_{n,k}(\tau_{k}^{+})=\alpha_{n}(\tau_{k})$. It follows from
(\ref{(eq): Boundary condition}) that $x_{n,i}^{\prime}(\tau_{k}^{+}%
)=x_{n,i}^{\prime}(\tau_{k}^{-})-h_{n}(\tau_{k})\tau_{k,i}^{\prime}$, $n=1,2$,
$i=1,2$. In case $(b)$, $x_{n}(\tau_{k}^{-})=0$, so that $f_{n,k-1}(\tau
_{k}^{-})=0$, and it is simple to verify that $x_{n,i}^{\prime}(\tau_{k}%
^{+})=x_{n,i}^{\prime}(\tau_{k}^{-})-\alpha_{n}(\tau_{k})\tau_{k,i}^{\prime}$,
$n=1,2$, $i=1,2$. Moreover, if the $k$th event corresponds to the $j$th
occurrence of a light switching event, we have $\tau_{k,i}^{\prime}%
=\sigma_{j,i}^{\prime}$ for some $j=1,2,\ldots$ Combining these results, we
get, for $n=1,2$ and $i=1,2$,%
\begin{equation}
x_{n,i}^{\prime}(\tau_{k}^{+})=x_{n,i}^{\prime}(\tau_{k}^{-})-\left\{
\begin{array}
[c]{c}%
h_{n}(\tau_{k})\sigma_{j,i}^{\prime}\\
\alpha_{n}(\tau_{k})\sigma_{j,i}^{\prime}%
\end{array}
\right.
\begin{array}
[c]{c}%
\text{if }x_{n}(\tau_{k})>0\\
\text{if }x_{n}(\tau_{k})=0
\end{array}
\label{(eq): G2R state derivative}%
\end{equation}
where $\sigma_{j,i}^{\prime}$ is given by (\ref{lemma1}) in Lemma 1 with
$\sigma_{j}=\tau_{k}$.

\begin{description}
\item[\textbf{(2)}] Event $R2G_{n}$
\end{description}

Once again, two cases must be considered: $(a)$ $R2G_{n}$ occurs at $\tau_{k}$
while road $n$ is undergoing an NEP; $(b)$ $R2G_{n}$ occurs at $\tau_{k}$
while road $n$ is undergoing an EP. In case $(a)$, the fact that road $n$ is
undergoing a RED cycle within a NEP at time $\tau_{k}^{-}$ means that
$f_{n,k-1}(\tau_{k}^{-})=\alpha_{n}(\tau_{k})$. Additionally, since road $n$
is undergoing a GREEN cycle at time $\tau_{k}^{+}$, we must have $f_{n,k}%
(\tau_{k}^{+})=\alpha_{n}(\tau_{k})-h_{n}(\tau_{k})$, and
(\ref{(eq): Boundary condition}) reduces to $x_{n,i}^{\prime}(\tau_{k}%
^{+})=x_{n,i}^{\prime}(\tau_{k}^{-})+h_{n}(\tau_{k}).\tau_{k,i}^{\prime}$,
$n=1,2$, $i=1,2$. In case $(b)$, the fact that road $n$ is empty while
undergoing a RED cycle at time $\tau_{k}^{-}$ implies that $f_{n,k-1}(\tau
_{k}^{-})=\alpha_{n}(\tau_{k})$ with $0<\alpha_{n}(\tau_{k})\leq h_{n}%
(\tau_{k})$, while $f_{n,k}(\tau_{k}^{+})=0$. Substituting these expressions
into (\ref{(eq): Boundary condition}) yields $x_{n,i}^{\prime}(\tau_{k}%
^{+})=x_{n,i}^{\prime}(\tau_{k}^{-})+\alpha_{n}(\tau_{k}).\tau_{k,i}^{\prime}%
$, $n=1,2$ and $i=1,2$. Combining these two cases, we get, for $n=1,2$ and
$i=1,2$,%
\begin{equation}
x_{n,i}^{\prime}(\tau_{k}^{+})=x_{n,i}^{\prime}(\tau_{k}^{-})+\left\{
\begin{array}
[c]{c}%
\alpha_{n}(\tau_{k})\sigma_{j,i}^{\prime}\\
\\
h_{n}(\tau_{k})\sigma_{j,i}^{\prime}%
\end{array}
\right.
\begin{array}
[c]{l}%
\text{if }x_{n}(\tau_{k})=0\text{ and }\\
0<\alpha_{n}(\tau_{k})\leq h_{n}(\tau_{k})\\
\text{otherwise}%
\end{array}
\label{(eq): R2G state derivative}%
\end{equation}
where again $\sigma_{j,i}^{\prime}$ is given by (\ref{lemma1}) in Lemma 1 with
$\sigma_{j}=\tau_{k}$.

\begin{description}
\item[\textbf{(3)}] Event $E_{n}$
\end{description}

This event corresponds to the end of an NEP on road $n$ and is induced by
$e_{5}$, which is an endogenous event at $\tau_{k}$ with $g_{k}=x_{n}=0$.
Since at time $\tau_{k}^{-}$ road $n$ is in an NEP, we must have
$f_{n,k-1}(\tau_{k}^{-})=\alpha_{n}(\tau_{k})-h_{n}(\tau_{k})$, and
(\ref{(eq): Endogenous event}) implies that $\tau_{k,i}^{\prime}%
=\frac{-x_{n,i}^{\prime}(\tau_{k}^{-})}{\alpha_{n}(\tau_{k})-h_{n}(\tau_{k})}%
$. Moreover, the fact that road $n$ is in an EP at time $\tau_{k}^{+}$ implies
that $f_{n,k}(\tau_{k}^{+})=0$, and (\ref{(eq): Boundary condition}) reduces
to $x_{n,i}^{\prime}(\tau_{k}^{+})=x_{n,i}^{\prime}(\tau_{k}^{-}%
)-x_{n,i}^{\prime}(\tau_{k}^{-})$ so that%
\begin{equation}
x_{n,i}^{\prime}(\tau_{k}^{+})=0\text{, \ \ \ }n=1,2\text{ and }%
i=1,2\label{(eq): Event E}%
\end{equation}

\begin{description}
\item[\textbf{(4)}] Event $S_{n}$
\end{description}

This event corresponds to the start of an NEP and can be induced by a
$G2R_{n}$, $e_{7}$ or $e_{6}$ event. These three cases are analyzed in what follows.

1. $S_{n}$\emph{ is induced by a }$G2R_{n}$\emph{ event}. Suppose that this
$G2R_{n}$ event initiated the $m$th NEP on road $n$. Therefore, during the
preceding EP, i.e. during the time interval $\left[  \eta_{n,m-1},\xi
_{n,m}\right)  $, we have $x_{n}(t)=0$ for $t\in\left[  \eta_{n,m-1},\xi
_{n,m}\right)  $, and, consequently, $x_{n,i}^{\prime}(t)=0$ for $t\in\left[
\eta_{n,m-1},\xi_{n,m}\right)  $ and $i=1,2$. As a result, $x_{n,i}^{\prime
}(\eta_{n,m-1}^{+})=x_{n,i}^{\prime}(\xi_{n,m}^{-})=0$, and since $\tau
_{k}=\xi_{n,m}$ it follows that $x_{n,i}^{\prime}(\tau_{k}^{-})=x_{n,i}%
^{\prime}(\xi_{n,m}^{-})=0$. Therefore, (\ref{(eq): G2R state derivative})
reduces to%
\begin{equation}
x_{n,i}^{\prime}(\tau_{k}^{+})=-\alpha_{n}(\tau_{k})\tau_{k,i}^{\prime
}\label{(eq): Event S case 1}%
\end{equation}
The value of $\tau_{k,i}^{\prime}$ above depends on the event inducing
$G2R_{n}$. If the $k$th event corresponds to the $j$th occurrence of a light
switching event, then $\tau_{k,i}^{\prime}=\sigma_{j,i}^{\prime}$ which is
obtained from (\ref{lemma1}). Note, however, that event $S_{n}$ cannot be
induced by $\gamma_{n}$ due to the fact that the occurrence of $\gamma_{n}$ is
conditioned upon road $n$ being in an NEP, which cannot be the case here. As a
result, the second case in (\ref{lemma1}) is excluded.

2. $S_{n}$\emph{ is induced by an }$e_{7}$\emph{ event}. Recall that $e_{7}$
corresponds to a switch from $\alpha_{n}(t)=0$ to $\alpha_{n}(t)>0$ while road
$n$ is undergoing a RED cycle, i.e. $z_{n}(t)=0$. Since this is an exogenous
event, $\tau_{k,i}^{\prime}=0$, $i=1,2$, and (\ref{(eq): Boundary condition})
reduces to $x_{n,i}^{\prime}(\tau_{k}^{+})=x_{n,i}^{\prime}(\tau_{k}^{-})$. We
know that $\tau_{k}$ corresponds to the time when the NEP starts at road $n$,
i.e. $\tau_{k}=\xi_{n,m}$, and we have shown that $x_{n,i}^{\prime}(\xi
_{n,m}^{-})=x_{n,i}^{\prime}(\eta_{n,m-1}^{+})=0$. It thus follows that
$x_{n,i}^{\prime}(\tau_{k}^{-})=x_{n,i}^{\prime}(\xi_{n,m}^{-})=0$, so that
$x_{n,i}^{\prime}(\tau_{k}^{+})=0$, $n=1,2$, $i=1,2$.

3. $S_{n}$\emph{ is induced by an }$e_{6}$\emph{ event}. Event $e_{6}$
corresponds to a switch from $\alpha_{n}(t)-h_{n}(t)\leq0$ to $\alpha
_{n}(t)-h_{n}(t)>0$ while road $n$ is undergoing a GREEN cycle, i.e.,
$z_{n}(t)>0$. Since this is an exogenous event, $\tau_{k,i}^{\prime}=0$,
$i=1,2$, and the subsequent analysis is similar to that of the previous case.
Therefore, $x_{n,i}^{\prime}(\tau_{k}^{+})=0$, $n=1,2$, $i=1,2$.

This completes the derivation of all state and event time derivatives required
to apply IPA to our TLC\ problem.

\subsection{Cost Derivatives}

Using the definition of $L(\mathbf{s})$ in
(\ref{(eq): Sample function simplified}) we can obtain the sample performance
derivatives $dL/ds_{i}$ as follows:%

\begin{gather*}
\frac{dL\left(  \mathbf{s}\right)  }{ds_{i}}=\frac{1}{T}\overset{2}%
{\underset{n=1}{\sum}}\underset{m=1}{\overset{M_{n}}{\sum}}%
{\displaystyle\int\limits_{\xi_{n.m}}^{\eta_{n,m}}}
w_{n}x_{n,i}^{\prime}\left(  t\right)  dt\\
+\frac{1}{T}\overset{2}{\underset{n=1}{\sum}}\underset{m=1}{\overset{M_{n}%
}{\sum}}\left[  w_{n}x_{n}\left(  \eta_{n,m}\right)  \frac{\partial\eta_{n,m}%
}{\partial s_{i}}-w_{n}x_{n}\left(  \xi_{n,m}\right)  \frac{\partial\xi_{n.m}%
}{\partial s_{i}}\right]
\end{gather*}
Note that $x_{n}\left(  \xi_{n,m}\right)  =x_{n}\left(  \eta_{n,m}\right)
=0$. Moreover, we have shown that $x_{n,i}^{\prime}(t)=x_{n,i}^{\prime}%
(\tau_{k}^{+})$, $t\in\left[  \tau_{k},\tau_{k+1}\right)  $, which implies
that we can decompose each NEP into time intervals of the form $\left[
\xi_{n,m},t_{n,m}^{1}\right)  ,\left[  t_{n,m}^{1},t_{n,m}^{2}\right)
,\ldots\lbrack t_{n,m}^{J_{n,m}},\eta_{n,m})$. Letting
\[
L_{n,m}(\mathbf{s})=\int\nolimits_{\xi_{n,m}}^{\eta_{n,m}}x_{n}(\mathbf{s}%
,t)dt
\]
we get%
\begin{align}
\frac{dL_{n,m}(\mathbf{s})}{ds_{i}} &  =x_{n,i}^{\prime}((\xi_{n,m})^{+}%
)\cdot(t_{n,m}^{1}-\xi_{n,m})\nonumber\\
&  +x_{n,i}^{\prime}((t_{n,m}^{J_{n,m}})^{+})\cdot(\eta_{n,m}-t_{n,m}%
^{J_{n,m}})\nonumber\\
&  +\sum\limits_{j=2}^{J_{n,m}}x_{n,i}^{\prime}((t_{n,m}^{j})^{+}%
)\cdot(t_{n,m}^{j}-t_{n,m}^{j-1})\label{IPAestimator}%
\end{align}
It is clear from (\ref{IPAestimator}) that computing the IPA estimator
requires knowledge of: $(i)$ the event times $\xi_{n,m}$, $\eta_{n,m}$, and
$t_{n,m}^{j}$, and $(ii)$ the value of the state derivatives $x_{n,i}^{\prime
}\left(  t\right)  $, whose expressions were derived in the previous section,
during each time interval. The quantities in $(i)$ are easily observed using
timers, and those in $(ii)$ ultimately depend on the values of the arrival and
departure rates $\alpha_{n}(t)$ and $h_{n}(t)$ at event times \emph{only},
which may be estimated through simple rate estimators. As a result, an
algorithm for updating the value of $dL\left(  \mathbf{s}\right)  /ds_{i}$
after each observed event is straightforward to implement. We also point out
that our IPA estimator is linear in the number of events in the SFM, not in
its states. Thus, our method can be readily extended to a network of intersections.

\section{Simulation Results}

With the intent of showing that performance improvements can be obtained when
IPA\ is used to control the queue content thresholds, two sets of simulations
were performed: one in which the thresholds were optimized considering a
priori fixed values of cycle lengths $\mathbf{\theta}=[\theta_{1,\min}%
,\theta_{1,\max},\theta_{2,\min},\theta_{2,\max}]$ for each road, and another
in which the cycle lengths and thresholds $\mathbf{s}=[s_{1},s_{2}]$ were
optimized sequentially. Thus, first, the IPA algorithm from \cite{Geng2013b}
was applied to determine optimal $\mathbf{\theta}$; then the values of $s_{1}$
and $s_{2}$ were optimized using the IPA algorithm described in this work.

In all our simulations, we assume that the vehicle arrival process is Poisson
with rate $\overline{\alpha}_{n}$, $n=1,2$, and approximate the departure rate
by a constant value $h_{n}(t)=H$ when road $n$ is non-empty, which amounts to
considering that the speed with which vehicles cross an intersection depends
only on the behavior of the vehicles themselves. We emphasize, however, that
our methodology applies independently of the distributions chosen to represent
the arrival and departure processes. We estimate the values of the arrival
rate at event times as $\alpha_{n}(\tau_{k})=N_{a}/t_{w}$, where $N_{a}$ is
the number of vehicle arrivals during a time window $t_{w}$ around $\tau_{k}$.
Simulations of the intersection modeled as a pure DES are thus run to generate
sample paths to which the IPA estimator is applied. We also make use of a
brute-force (BF) approach to generate a cost surface along which the
convergence of the IPA-driven optimization algorithm is depicted. The BF
method consists of discretizing the values of $s_{i}$ and generating 10 sample
paths for each pair of discretized threshold values $(s_{1},s_{2}%
)=(1,1),(1,2),\ldots,(2,1),\ldots$, from which the average total cost can then
be obtained. In all results reported here, we set $H=1$, $w_{n}=1$, $n=1,2$,
and measure the sample path length in terms of the number of observed light
switches, which we choose to be $N=5000$.

\begin{figure}
\begin{center}%
$\includegraphics[width=8.4cm]{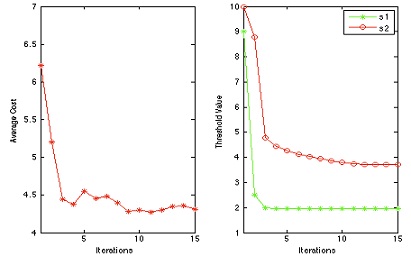}$\caption{Sample cost and
parameter trajectories for $1/\bar{\alpha}=[2,6]$, $\theta=[10,30,10,30]$,
and $s_0=[9,10]$}\label{(fig): Convergence ARB 2}%
\end{center}
\end{figure}

\begin{figure}
\begin{center}%
$\includegraphics[width=8.4cm, height=8.4cm]{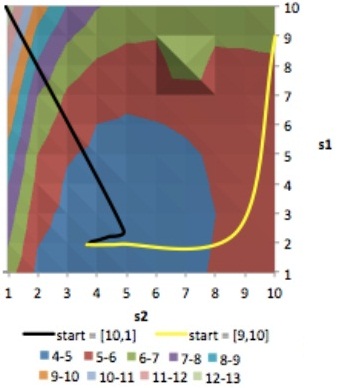}$%
\caption{Cost surface and
convergence trajectories for \emph{Scenario A} (note: the color scale refers
to the cost values)}\label{(fig): Scenario A}%
\end{center}
\end{figure}

In our first set of simulations, the GREEN light cycles are fixed and equal on
both roads by setting $\theta_{n,\min}=10$ $\sec$ and $\theta_{n,\max}=30$
$\sec$, $n=1,2$. Two scenarios are considered: $Scenario$ $A$, in which road
$1$ exhibits high traffic intensity while road $2$ exhibits low traffic
intensity: $1/\overline{\alpha}_{1}=2$ and $1/\overline{\alpha}_{2}=6$;
$Scenario$ $B$, in which both roads exhibit high but unequal traffic
intensity: $1/\overline{\alpha}_{1}=2$ and $1/\overline{\alpha}_{2}=3$. We
further consider two different initial threshold configurations for each
scenario. Table \ref{(tb): Results arbitrary theta} shows the optimal
threshold values determined by both the BF method and the IPA-driven
optimization algorithm, along with the cost reduction achieved by the latter
(denoted as $R$ and computed as a percentage of the initial cost). Sample
convergence plots of the cost $J$ and thresholds $\mathbf{s}$ are presented in
Fig. \ref{(fig): Convergence ARB 2}, while the cost surface referring to
$Scenario$ $A$, along with curves (black and yellow) that represent the
trajectories corresponding to each initial configuration, is shown in Fig.
\ref{(fig): Scenario A}. Visual inspection of Fig. \ref{(fig): Scenario A}
reveals that both trajectories converge to the same optimal point, namely
$s_{IPA}^{\ast}=[1.9,3.7]$, as presented in Table
\ref{(tb): Results arbitrary theta}.

\begin{table}[bh]
\caption{Optimization results for system with a priori fixed cycle lengths}%
\label{(tb): Results arbitrary theta}
\begin{center}%
\begin{tabular}
[c]{|c|c|c|c|c|c|c|c|}\hline
& \multicolumn{2}{|c|}{Initial Point} & \multicolumn{3}{|c|}{IPA} &
\multicolumn{2}{|c|}{BF}\\\hline
$1/\overline{\alpha}$ & $s_{0}$ & $J_{0}$ & $s_{IPA}^{\ast}$ & $J_{IPA}^{\ast
}$ & $R$ & $s_{BF}^{\ast}$ & $J_{BF}^{\ast}$\\\hline
$\left[  2,6\right]  $ & $\left[  10,1\right]  $ & $12.8$ & $\left[
1.9,3.7\right]  $ & $4.3$ & $66$ & $\left[  1,4\right]  $ & $4.4$\\\hline
$\left[  2,6\right]  $ & $\left[  9,10\right]  $ & $6.2$ & $\left[
1.9,3.7\right]  $ & $4.3$ & $31$ & $\left[  1,4\right]  $ & $4.4$\\\hline
$\left[  2,3\right]  $ & $\left[  15,3\right]  $ & $18.9$ & $\left[
4.6,5.1\right]  $ & $7.9$ & $58$ & $\left[  5,6\right]  $ & $8.8$\\\hline
$\left[  2,3\right]  $ & $\left[  15,15\right]  $ & $13.1$ & $\left[
4.6,5.1\right]  $ & $7.9$ & $40$ & $\left[  5,6\right]  $ & $8.8$\\\hline
\end{tabular}
\end{center}
\end{table}%

%

In our second set of simulations, we perform a sequential optimization of the
cycle lengths and threshold values. We make use of the optimal light cycle
lengths obtained through IPA (denoted by $\mathbf{\theta}_{IPA}^{\ast}=\left[
\theta_{1,\min}^{\ast},\theta_{1,\max}^{\ast},\theta_{2,\min}^{\ast}%
,\theta_{2,\max}^{\ast}\right]  $) by \cite{Geng2013b}, and subsequently apply
the IPA estimator derived in this paper to optimize the queue content
thresholds. The optimal light cycle lengths obtained by \cite{Geng2013b} for
fixed and predetermined threshold values of $\mathbf{s}=\left[  8,8\right]  $
were $\mathbf{\theta}_{IPA}^{\ast}=\left[  10.2,19.3,10.1,16.3\right]  $ for
$1/\overline{\alpha}=[2,3]$ and $\mathbf{\theta}_{IPA}^{\ast}=\left[
10.1,20.1,10.6,11.9\right]  $ for $1/\overline{\alpha}=[1.7,3]$. A comparison
between IPA and BF results, including a quantitative assessment of the
\emph{additional} cost reduction achieved (computed as a percentage of the
initial cost and labeled $R$) is shown in Table
\ref{(tb): Results sequential optimization}.

\begin{table}[bh]
\caption{Optimization results for system with optimal cycle lengths}%
\label{(tb): Results sequential optimization}
\begin{center}%
\begin{tabular}
[c]{|c|c|c|c|c|c|}\hline
& \multicolumn{3}{|c|}{IPA} & \multicolumn{2}{|c|}{BF}\\\hline
$1/\overline{\alpha}$ & $s_{IPA}^{\ast}$ & $J_{IPA}^{\ast}$ & $R$ $(\%)$ &
$s_{BF}^{\ast}$ & $J_{BF}^{\ast}$\\\hline
$\left[  2,3\right]  $ & $\left[  2.8,4.3\right]  $ & $7.1$ & $15$ & $\left[
2,5\right]  $ & $7.2$\\\hline
$\left[  1.7,3\right]  $ & $\left[  4.8,6.1\right]  $ & $14.9$ & $11$ &
$\left[  3,8\right]  $ & $15.7$\\\hline
\end{tabular}
\end{center}
\end{table}

In order to further illustrate the advantage of quasi-dynamically controlling
the light cycle lengths and threshold values over a static IPA approach to the
TLC problem, we include a comparison of the results generated by our
methodology with those obtained when static control (as described by
\cite{Geng2012}) is applied to determine the optimal cycle lengths
$\mathbf{\theta}_{static}^{\ast}$. The static controller defined by
\cite{Geng2012} adjusts the green light cycles subject to some lower and upper
bounds and determines $\mathbf{\theta}_{static}^{\ast}=\left[  \theta
_{1}^{\ast},\theta_{2}^{\ast}\right]  $, where $\theta_{1}^{\ast}$
($\theta_{2}^{\ast}$, respectively) is the green cycle length which should be
allotted to road $1$ (road $2$, respectively) in order to minimize the average
queue content on both roads. Table \ref{(tb): Results comparison} summarizes
the results obtained by each of the IPA approaches considered in this work:
\textit{Method 1}, in which a static controller is used to adjust the light
cycles (results were obtained by using the same setting as in our second set
of quasi-dynamic simulations and constraining $\theta\in\left[  10,40\right]
$); \textit{Method 2}, in which only the light cycles are controlled
quasi-dynamically (i.e. fixed and predetermined queue content thresholds are
incorporated into the system model); \textit{Method 3}, in which a sequential
quasi-dynamic optimization of light cycle lengths and threshold values is
performed in between two adjustment points. The columns labeled $R_{i}$,
$i=2,3$, present the cost reduction achieved by the quasi-dynamic methods with
respect to the static approach, i.e. $R_{i}=\frac{J_{1}^{\ast}-J_{i}^{\ast}%
}{J_{1}^{\ast}}\ast100$.

\begin{table}[bh]
\caption{Comparison between three IPA-based approaches to the TLC problem}%
\label{(tb): Results comparison}
\begin{center}%
\begin{tabular}
[c]{|c|c|c|c|c|c|}\hline
& Method 1 & \multicolumn{2}{|c|}{Method 2} & \multicolumn{2}{|c|}{Method
3}\\\hline
$1/\overline{\alpha}$ & $J_{1}^{\ast}$ & $J_{2}^{\ast}$ & $R_{2}$ $(\%)$ &
$J_{3}^{\ast}$ & $R_{3}$ $(\%)$\\\hline
$\left[  2,3\right]  $ & $14.4$ & $8.4$ & $42$ & $7.1$ & $51$\\\hline
$\left[  1.7,3\right]  $ & $23.9$ & $16.7$ & $30$ & $14.9$ & $38$\\\hline
\end{tabular}
\end{center}
\end{table}

\section{Conclusion}

We have modeled a single traffic light intersection as an SFM and formulated
the corresponding TLC problem within a quasi-dynamic control setting to which
IPA techniques were applied in order to derive gradient estimates of a cost
metric with respect to controllable queue content threshold values. By
subsequently incorporating these estimators into a gradient-based optimization
algorithm, numerical results were obtained which substantiate our claims that:
$(i)$ a considerable reduction in the mean queue content of both roads can be
achieved by quasi-dynamically controlling the thresholds in systems with
non-optimal cycle lengths; $(ii)$ determining optimal threshold values allows
for additional improvements to the performance of systems running under
optimal cycle lengths. Such results lead us to believe that a method in which
the light cycle lengths and queue content thresholds are controlled
\emph{simultaneously} is likely to provide improved solutions to the
TLC\ problem. Our ongoing research is, therefore, focused on deriving an
IPA-based optimization algorithm that incorporates all such controllable
parameters and ultimately determines the optimal light cycle length/threshold
configuration capable of minimizing traffic build-up at a given traffic light
intersection. Future work includes applying IPA to an intersection with more
complicated traffic flow (e.g. allowing for left- and right-turns),
incorporating acceleration/deceleration due to light switches into the flow
model, as well as extending our methodology to multiple intersections.

\end{document}